\documentclass[12pt]{spieman}  
\pdfoutput = 1
\usepackage{amsmath,amsfonts,amssymb,lineno}
\usepackage{algorithm,setspace}
\usepackage{algpseudocode}
\usepackage{graphicx}
\usepackage{setspace}
\usepackage{tocloft}

\usepackage{color}
\definecolor{gray}{rgb}{0.6,0.6,0.6}
\definecolor{red}{rgb}{0.85,0,0}
\definecolor{green}{rgb}{0,0.85,0}
\definecolor{blue}{rgb}{0,0,0.85}
\definecolor{beige}{rgb}{0.92,0.87,0.78}

\DeclareMathOperator*{\argmin}{arg\,min}

\title{Improving depth-resolution, in-plane contrast, and reducing non-uniformity artifacts for wide-angle DBT}

\author[a,*]{Emil Y. Sidky}
\author[b]{Adrian A. Sanchez}
\author[a]{John Paul Phillips}
\author[a]{Zheng Zhang}
\author[a]{Dan Xia}
\author[a]{Ingrid S. Reiser}
\author[a]{Xiaochuan Pan}
\affil[a]{The University of Chicago, Department of Radiology MC-2026, 5841 S. Maryland Ave., Chicago, IL 60615}
\affil[b]{ Washington University , Mallinckrodt Institute of Radiology, Campus
Box 8131‑19‑01, 510 South Kingshighway Boulevard,
Hazelwood, St. Louis, MO 63110}

\cftpagenumbersoff{figure}
\cftpagenumbersoff{table} 
\begin{document} 
\maketitle

\begin{abstract}
\\
\textbf{Purpose: } This work aims to develop an image reconstruction algorithm for wide-angle
digital breast tomosynthesis (DBT) that has improved depth resolution and in-plane contrast
while reducing non-uniformity artifacts.\\
\textbf{Approach: } 
The image reconstruction algorithm is an extension of our prior work on sparsity-regularized
iterative image reconstruction. The algorithm is performed in two stages as explained in a
prior work. The first stage consists of a low-resolution reconstruction that exploits sparsity
for quantitative accuracy. In this work, this first stage is augmented with a formulation that
includes the estimation of a ``background'' image, which absorbs low-frequency artifacts
that cause image non-uniformity.\\
\textbf{Results: } The new algorithm is demonstrated on three patient cases, each with
different features of clinical interest, and for which the data are
acquired on a wide-angle DBT system.\\
\textbf{Conclusion: } The results on the shown cases indicate that the algorithm design goals have
been met, but additional empirical results and task-based assessment are needed to strengthen
this conclusion.
\end{abstract}


{\noindent \footnotesize\textbf{*}Emil Y. Sidky,  \linkable{sidky@uchicago.edu} }

\begin{spacing}{2}   

\section{Introduction}
\label{sec:intro}

The impact of Professor Barrett's work on our research is profound and multi-faceted and it would be impossible
to detail every instance where his influence was felt.
Here, in this special issue article, we focus specifically on our recent line of investigation
into sparsity-exploiting image reconstruction for digital breast tomosynthesis (DBT).
In this project, we have mainly made use of Chapter 7
``Deterministic Description of Imaging Systems'' in the {\it Foundations of Image Science} (FIS)\cite{barrett2003foundations}, 
coauthored by Professor Barrett and Dr. Kyle Myers.
Understanding the details of this chapter is crucial to anyone who is interested in developing
image reconstruction algorithms, because image reconstruction is essentially an inversion of an imaging system model
and thus a clear understanding of the types of image system models is necessary.

The organization of the paper is as follows: Sec.~\ref{sec:imaging} explains how the taxonomy of imaging system models
is utilized to formulate our studies on iterative image reconstruction for limited-data systems such as DBT.
Sec.~\ref{sec:coupled} then develops an extension to our previous DBT work that addresses depth resolution,
in-plane contrast, and non-uniformity artifacts; Sec.~\ref{sec:highres} presents methods to achieve standard DBT in-plane
resolution; Sec.~\ref{sec:setup} describes the DBT scanner configuration along
with the algorithm parameter settings,
and Sec.~\ref{sec:results} demonstrates the algorithm on three patient
cases acquired on a wide-angle DBT system.

\section{Imaging system models and sparsity-exploiting image reconstruction in digital breast tomosynthesis}
\label{sec:imaging}

Chapter~7 in FIS describes
three imaging system models: continuous-to-continuous (CC), continuous-to-discrete (CD), and discrete-to-discrete (DD).
For X-ray-based tomography,
the divergent-ray and Radon transforms are instances of CC models, where the object model is a continuous
function of spatial variables and the measurement space is imagined to be a continuous function of X-ray source angle
and detector receptor location.
The CD model, however, is a more accurate description for any form of digital imaging because digital
detectors consist of discrete arrays of pixels. Finally, the DD model must be contended with because image processing
or reconstruction algorithms are implemented on computers, which necessarily discretize images so that they can be stored
in memory and manipulated by the various programmed processing units. In addition, fully utilizing Dr. Barrett's taxonomy
of imaging models and
his emphasis on task-based image quality assessment is reflected throughout all of our work on imaging.

For analytic image reconstruction, the filtered back-projection (FBP) formulas are derived to be the inverse of the
CC divergent-ray or Radon transforms. Because the measured data are discrete, some form of interpolation is
implicit when creating the input sinogram to an FBP algorithm. Because interpolation may not perfectly reconstruct
the continuous sinogram, the output image values will have error even for ideal, noiseless sinogram measurements.
Furthermore, the FBP algorithm can only be executed on a computer for a finite spatial grid. Thus,
interpolation is again needed if the display pixel grid does not match the grid of points, which has been reconstructed
by FBP. For iterative image reconstruction, the DD model plays a larger role because the scanned object is represented
by a finite set of expansion functions, which is most commonly taken to be a grid of voxels. Thus, another source
of error is introduced; namely, the divergent-ray transform of a grid of voxels likely does not match the divergent-ray
transform of a continuous object function for any set of voxel coefficients. Thus discretization error is an issue
with which all iterative algorithms must contend.

Turning to sparsity-exploiting image reconstruction in X-ray tomography,
we write the canonical optimization problem for compressive
sensing (CS) \cite{candes2006robust}
\begin{equation}
\label{CS}
\argmin_f \| \nabla f \|_1 \;\; \text{such that} \;\; Xf=g,
\end{equation}
where $f$ represents the image; $g$ is the sinogram; $X$ is the X-ray, or divergent-ray, transform; and
the objective function is the image total variation (TV). This optimization is designed to recover the image
$f$ under the assumption that $f$ has a sparse gradient and the sinogram data is ideal; namely the data are generated
from the true underlying object $f_0$, i.e. $g = Xf_0$. This CS optimization becomes interesting when the
measured data coverage is less than what would be needed for FBP image reconstruction; i.e. the equation
$Xf = g$ does not have a unique solution. As of yet, we have not
specified the type of imaging model; this optimization can be expressed in CC, CD, and DD forms.
It turns out that the vast
majority of the work on CS does in fact assume a DD imaging system model;
the conditions for recovery of $f$, discussed in the literature, mainly 
focus on the number of samples in $g$ being some factor
larger than the number of non-zeros in the image gradient. Such conditions only make sense with a DD formulation.
Accordingly, $\nabla$ becomes a matrix encoding a finite differencing approximation to the continuous gradient;
$X$ is the discretized X-ray transform; $g$ is a set of sinogram samples; and $f$ is a vector of image basis coefficients,
most commonly voxel values.

For X-ray tomography, the CS optimization problem was originally formulated in CT to address view-angle undersampling,
limited angular-range scanning, and cone-beam artifacts from a circular scan,  \cite{sidky2006accurate,sidky2008image}
but as explained in Pan {\it et al.} \cite{pan2009commercial} CS can address a number of under-sampling issues such as projection truncation,
bad detector pixels, and exterior/interior tomography. An important point to understand about all of these sampling
challenges is that their degree of difficulty can be tuned for DD imaging system models. Because much CS work
focuses on data sampling, mentioning the image basis expansion in passing or perhaps not at all, it leaves
the impression that the sampling challenge is being addressed for CC or CD imaging system models and that a continuous
image function -- or a close approximation thereof -- is being recovered.
The reality is, however, that DD imaging system models are being employed, and only a finite
set of image basis coefficients are being recovered. With this being the case any sampling challenge can be
overcome simply by reducing the size of the basis expansion set. The price to be paid is that discretization
error increases when performing image reconstruction on real data or simulation data generated with a continuous image model.

These issues surrounding DD modeling are central to our work on the limited angular-range image reconstruction problem that arises in DBT.
Our first steps in accurately solving a DD model for DBT were reported in Zhang {\it et al.} \cite{zhang2021directional}, where
the optimization problem of interest was a variant of TV-based image reconstruction that used directional TV (DTV) constraints.
Following up on this work, a thorough investigation on wide-angle DBT image reconstruction with DTV penalties
was presented in Sidky {\it et al.} \cite{sidky2025accurate}. The wide-angle system studied was that of the Siemens Mammomat
scanner that has 85 micron detector pixel resolution and scans the breast with 25 exposures
covering a 50 degree scan arc. The first
study in that paper directly addressed the issue of ``tuning'' the DD model in order to get accurate image reconstruction
for wide-angle DBT. Based on those studies, it was determined that relatively accurate and stable image reconstruction could be
performed in a volume composed of cubic voxels with a width four times the width of the detector pixels. The remainder
of that work investigated the impact of discretization error when the data are generated from a CD breast phantom model,
and additionally the DTV algorithm was applied to DBT image reconstruction for a physical phantom obtained on a Mammomat
scanner.

Central to all of these studies is the need for accurate solvers for non-smooth convex optimization problems.
Many of the DD imaging system models of interest do not have analytic solutions available and they involve convex
optimization problems with non-smooth terms such as the DTV norm or hard constraints such as non-negativity.
The only way to provide evidence that the DD imaging system model of interest can be solved is to gather empirical
results with computer-simulated phantoms, where the ground truth is known exactly. It must then be demonstrated
that the phantom can be recovered perfectly from ideal noiseless data generated with this phantom. For this reason,
we have been consistently adapting and refining the primal-dual hybrid gradient (PDHG) algorithm \cite{chambolle2011first}
for X-ray tomographic image reconstruction, \cite{sidky2012convex} where the latest development is focused
on efficient implementation for multi-term optimization problems, which are specifically needed for DD imaging system
models relevant to DBT. \cite{sidky2026efficient}

The imaging system model design in this latest work is heavily influenced by the ideas discussed in FIS.
The overall design is to split the reconstruction problem into two stages: a low-resolution sparsity-exploiting
reconstruction using a least-squares
(LSQ) data fidelity with DTV regularization (lowres-LSQ-DTV), and high resolution
Tikhonov-regularized least-squares (highres-LSQ-Tik).  The design goal of lowres-LSQ-DTV is to minimize image
root-mean-square-error (RMSE) at low resolution. The low-resolution aspect is needed because "inverse crime"
inversion of the DD imaging model tells us that accurate image reconstruction becomes possible for cubic
voxels four times larger than the detector pixels for wide-angle DBT as demonstrated in Refs.~\citenum{sidky2025accurate,sidky2026efficient}.
To avoid artifacts at this low-resolution both the object model and input data make use of gaussian smoothing
that helps to avoid a jagged image representation and to ensure that the data fed into lowres-LSQ-DTV do not
contain frequency information that is not representable with the low-resolution image grid; the need for this
data filtering becomes clear from the various forms of discretization error explained in Chapter 7 of FIS.
The use of large voxels and the fact that lowres-LSQ-DTV is optimized on image RMSE means that small, yet
important, features in the image could be lost due to over-regularization -- as demonstrated
by a signal detection study using the channelized Hotelling observer machinery that was so well-explained in FIS \cite{sidky2021signal}.
To bring in the necessary fine details into the final image, the lowres-LSQ-DTV output image is upsampled to $f_\text{DTV}$ and
this image initializes the highres-LSQ-Tik algorithm, where the Tikhonov penalty terms also ensures proximity to $f_\text{DTV}$.
Because lowres-LSQ-DTV provides a low RMSE initial image, the parameters of highres-LSQ-Tik can then be objectively determined
based solely on clinical tasks of interest to DBT such as the detection of spiculations or microcalcifications using the machinery
developed in FIS.

\section{Coupled lowres-LSQ-DTV image reconstruction for estimating background artifacts}
\label{sec:coupled}

In this section, we present another advance of our DBT image reconstruction framework that deals with low-frequency
background artifacts. Such artifacts enter into the raw reconstruction for a variety of reasons: imperfect knowledge
of the X-ray illumination in each projection, X-ray scatter, and blurred high-density regions above or below the viewing plane. 
These low-frequency artifacts complicate display of the in-plane slices because they limit the ability to set
the gray scale window for maximizing contrast of features of interest.
In our previous work,
Ref.~\citenum{sidky2025accurate} reconstructed in-plane slices were post-processed to remove low-frequency
artifacts by fitting a low-order 2D polynomial to each slice and dividing the raw slice image by this fit.
While effective for image presentation, it is not a 3D solution and it loses quantitative information.

Here, we augment lowres-LSQ-DTV, developed in Ref.~\citenum{sidky2026efficient}, to coupled lowres-LSQ-DTV where
three images are reconstructed simultaneously. The raw reconstructed image is denoted by $f_1$; the ``background''
image, which will primarly contain the low-frequency artifacts, is denoted by $f_2$; and the difference image
is $f_3 = f_1 - f_2$. The optimization problem for coupled lowres-LSQ-DTV is
\begin{align}
\argmin_{f_1,f_2,f_3} \Bigl\{
 &  \alpha_x \| \partial_x f_1 \|_1 + \alpha_y \| \partial_y f_1 \|_1  +
  \alpha_a \| \partial_a f_1 \|_1  + \alpha_b \| \partial_b f_1 \|_1 + \alpha_1 \| f_1 \|_1 \notag \\
+&  \alpha_x \| \partial_x f_2 \|_1 + \alpha_y \| \partial_y f_2 \|_1  +
  \alpha_a \| \partial_a f_2 \|_1  + \alpha_b \| \partial_b f_2 \|_1 + \alpha_1 \| f_2 \|_1 \notag \\
+&   \alpha_3 \| f_3 \|_1 \;\;\; \text{ such that } \;\;\; f_1 - f_2 = f_3, \notag \\
 &  \left\|R[c] \left(XG[d_1]f_1-g\right) \right\|_2 \le \epsilon_1 \cdot \sqrt{N_\text{data}}, \text{ and } \notag\\
 & \left\|R[c] \left(XG[d_2]f_2-G_\text{det}[d_d]g\right) \right\|_2 \le \epsilon_2 \cdot \sqrt{N_\text{data}}  \Bigr\}, \label{DTVcopt3D}
\end{align}
where $g$ represents the DBT projection data and $N_\text{data}$ is the size of this dataset; $X$ is the DBT X-ray transform
that computes a projection data estimate based on an image estimate $f$; the operator $G[d]$ performs 3D gaussian blurring with
width $d = (d_x,d_y,d_z)$ and likewise $G_\text{det}[d]$ performs 2D gaussian blurring on each projections with width $d=(d_u,d_v)$;
the operator $R[c]$ performs projection filtering with the square root of the ramp filter with frequency cutoff parameter $c$;
the parameters $\epsilon_1$ and $\epsilon_2$ are RMSE constraint parameters on the projection estimate for images $f_1$ and $f_2$;
the objective function consists of a number of sparsity-inducing $\ell_1$-norm terms; the partial derivative symbol represents
the finite differencing approximation to the derivative in the direction indicated in the subscript, where $a$ and $b$ represent
oblique directions corresponding to vectors perpendicular to first and last X-ray source position in the $xz$-plane, see
Sec.~\ref{app:pdhg};
the parameters $\alpha_x$, $\alpha_y$, $\alpha_a$, $\alpha_b$, and $\alpha_1$ specify the relative
strength of each sparsity-inducing term by normalizing the sum of these parameters to 1; note that the same parameters
are used for $f_1$ and $f_2$; and finally, the parameter $\alpha_3$ sets the strength on the $\ell_1$ penalty for the difference
image $f_3$. The PDHG algorithm for this optimization is derived by mapping it onto the generic PDHG form; Sec~\ref{app:pdhg}
shows how this mapping is done and it specifies the step-size parameterization.

A large part of the image model design leading to Eq.~(\ref{DTVcopt3D}) is explained in our previous paper Ref.~\citenum{sidky2026efficient}.
The new aspect is setting up the model for simultaneous optimization of the raw image $f_1$ and the ``background'' image $f_2$ and coupling
these two images through the $\ell_1$-norm penalty on $f_3$.
Use of two data fidelity terms with low- and high-frequency bands has the potential to improve algorithm convergence\cite{iyadomi2026}
in addition to the estimation of $f_1$ and $f_2$ images.
The use of the $\ell_1$-norm penalty on $f_3$ is to encourage small differences between $f_1$
and $f_2$ to continue to decrease, and more importantly, the $\ell_1$-norm does not strongly penalize large differences. We are expecting
that $f_1$ will deviate from $f_2$ for in-focus features of the breast. For the data fidelity term involving $f_2$, the intention for
introducing the gaussian blur $G[d_2]$ is to smooth out in-plane structures; accordingly the $x$ and $y$ components of the width vector $d_2$
should be large, on the order of centimeters, while the depth or $z$ component of $d_2$ should be on the order of the slice thickness.
The corresponding blur for projection data, $G_\text{det}[d_d]$, should be slightly larger than the in-plane blurring components of $d_2$
in order to suppress frequency components that would pick up real structures in the breast and so that filtered data can be consistent
with the designed imaging model for $f_2$. For the data fidelity term involving $f_1$, the purpose of the gaussian blur $G[d_1]$ is to
smooth the blockiness introduce by the use of large voxels before X-ray projection.

\section{Achieving standard in-plane DBT resolution}
\label{sec:highres}
The coupled lowres-LSQ-DTV is designed to greatly improve quantitative accuracy at low-resolution, which has the side benefit
of improving depth-resolution. This by itself may have clinical use for dense tissue mapping and refined dose models, but for
standard DBT image viewing it is necessary to achieve in-plane resolution close to that of the detector. For improving
in-plane resolution, two methods are investigated.

First, we follow the strategy
of Ref.~\citenum{sidky2026efficient} and use $f_1$ to initialize a high-resolution Tikhonov-regularized least-squares (DTVinit-LSQ-Tik) image
reconstruction
\begin{equation}
\label{LSQTikHR}
\argmin_h \left\{ \tfrac{1}{2}\|R[c](XG[d]h-g)\|^2_2 + \alpha_\text{Tik} \tfrac{1}{2} \| G[d]h - h_0 \|^2_2  \right\},
\end{equation}
where the high-resolution image and upsampled initial image are represented by $h$ and $h_0$, respectively.
Ref.~\citenum{sidky2026efficient} explains
how to form the initial image $h_0$ from the low-resolution image $f_1$; and the gaussian blurring width for DTVinit-LSQ-Tik
is chosen to be the same as the voxel dimensions. The algorithm consists of steepest descent steps on the LSQ-Tik objective function.
The parameter $\alpha_\text{Tik}$ controls the image regularization.
In order to form the final image for display, the background image $f_2$ is upsampled to the high-resolution background image $h_2$.
The support image, $h_\text{sup}$, indicating voxels inside and outside of the breast with 1 and 0, respectively, is found by
thresholding $h_2$ at a gray level of $0.1$ cm$^{-1}$. A mean background gray value, $b$, is found by averaging $h_2$ in the breast support,
i.e. where $h_\text{sup} = 1$. The image for display is then formed by
\begin{equation}
\label{disp1}
h_\text{disp-LSQ} = G_z[d_z] \left((h - h_2) + b \cdot h_\text{sup} \right),
\end{equation}
where additional depth blurring with a one-dimensional gaussian kernel is introduced
for additional regularization.

Second, the coupled lowres-LSQ-DTV result is reprojected and subtracted from the original DBT projection data. This difference
is reconstructed with FBP to form an image that contains only the high-resolution detail, which can subsequently be combined
with the coupled lowres-LSQ-DTV result. The processing for the high-resolution DTV-FBP algorithm starts with forming the high-resolution
difference image
\begin{equation}
\label{DTVFBP}
h_\text{diff} = \text{FBP}(g - X G[d_1] h_1),
\end{equation}
where the high-resolution image $h_1$ is obtained by upsampling
the raw coupled lowres-LSQ-DTV image, $f_1$; the Gaussian blurring $G[d_1]$ uses a width $d_1$ matched to the value employed for 
Eq.~(\ref{DTVcopt3D}). In our implementation of FBP, a Hanning filter is used with cutoff set to the Nyquist frequency.
For this method, the display image is formed slightly differently than what is described in Eq.~(\ref{disp1}).
\begin{equation}
\label{disp2}
h_\text{disp-FBP} = G_z[d_z] \left((h_1 + \eta \cdot  h_\text{diff} - h_2) + b \cdot h_\text{sup} \right),
\end{equation}
where the combination parameter $\eta$ needs to be tuned. The useful scale for $\eta$ can be determined by reprojecting
$h_\text{diff}$ and normalizing to the magnitude of the projection data difference, $g - X G[d_1] h_1$.

\section{Scanner setup and algorithm parameter settings}
\label{sec:setup}
The algorithm parameter settings depend on the configuration of the DBT scanner, which can vary dramatically
from vendor to vendor. For this work, patient data from the Siemens Mammomat B.brilliant scanner is reconstructed.
The scan configuration
involves 25 X-ray projections covering an arc of 50 degrees, which is currently the widest angle DBT on the market. The X-ray detector
is 3584$\times$2816 pixels with each pixel being 85 microns squared. The long-axis of the detector $u$ is parallel to the direction
of the scan and the $x$-axis of the image volume, while the short-axis $v$ is parallel to the $y$-axis of the image volume.
To process the transmission data to projection data, the air scan fluence, $I_0(u,v)$, is estimated from a strip on the detector
well beyond the breast projection.
The $u$ dependent fluence,  $\bar{I}_0(u)$, is estimated by averaging $I_0(u,v)$ over $v$ only within this strip.
The estimated projection data, $g(u,v)$, are obtained by standard negative logarithm processing;
i.e. $g(u,v) = - \ln(I(u,v)/\bar{I}_0(u))$.

The coupled lowres-LSQ-DTV image reconstruction is carried out on a low-resolution
voxel grid with voxels of size 340x340x680 microns cubed, and the projection data are
then down-sampled to a pixel grid with pixels of size
340x340 microns squared. The high-resolution LSQ-Tik and DTV-FBP algorithms use
voxels of size 85x85x340 microns cubed, and the native
85x85 microns squared detector resolution for the projections is maintained.

The regularization parameter settings with which coupled lowres-LSQ-DTV is solved are
\begin{equation*}
\alpha_y = \alpha_a = \alpha_b = \alpha_1 =1/9; \; \; \alpha_x = 5/9,
\end{equation*}
which are the same values determined by DBT simulations on the same wide-angle DBT system,
performed in Ref.~\citenum{sidky2026efficient},
and the coupling parameter for the difference image $f_3$ is
set to
\begin{equation*}
\alpha_3 = 0.2,
\end{equation*}
which provides a good visual fit to the image background for the three patient cases shown in Sec.~\ref{sec:results}.
The blur parameters $d_1$, $d_2$, and $d_d$ from the data discrepancy constraints of Eq.~(\ref{DTVcopt3D}) are set
to
\begin{equation*}
d_1 = (340,340,680) \; \mu m, \; \; d_2 = (3400,3400,680) \; \mu m, \; \;\text{and } d_d=(3740,3740) \; \mu m,
\end{equation*}
for all three patient cases. That the same parameters can be applied to each of these cases is an indication
that the imaging is not so sensitive to these parameter settings.

\begin{table}[h!]
\centering
\begin{tabular}{c|| c | c | c }
case & $\epsilon_1$ & $\epsilon_2$ & $d_z$ \\ [0.5ex]
\hline\hline
mass  & 0.015 & 0.015 & 3400 $\mu$m  \\
\hline
distortion  & 0.012 & 0.01 & 2040 $\mu$m  \\
\hline
calcs & 0.02  & 0.0125  & 1360 $\mu$m

\end{tabular}
\caption{Parameter settings for the coupled lowres-LSQ-DTV optimization and image display depth blurring for each of
the patient cases. The cases presented in order contain a solid mass (mass), an architectural distortion (distortion),
and a calcification cluster (calcs).}
\label{tab:parms}
\end{table}

\begin{figure}
\begin{center}
\begin{tabular}{c}
\includegraphics[width=0.9\columnwidth]{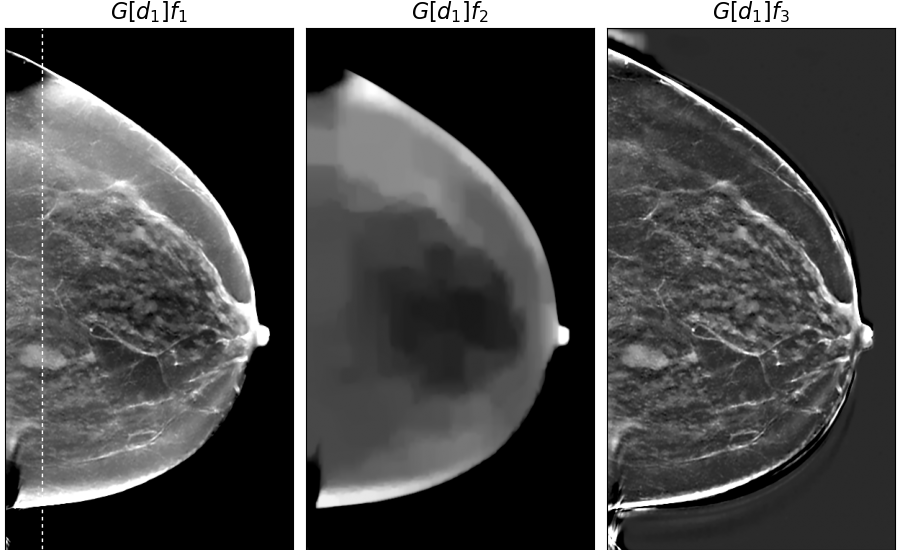} \\
\includegraphics[width=0.5\columnwidth]{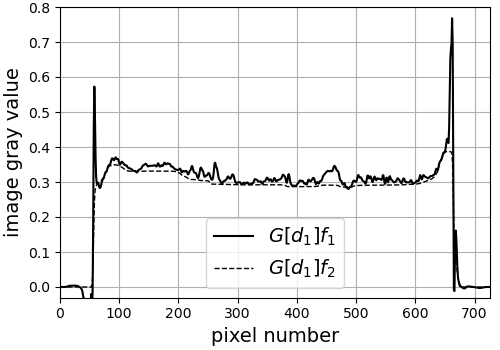}
\end{tabular}
\end{center}
\caption[Low-resolution coupled lowres-LSQ-DTV mass images.]{ \label{fig:DTV-LSQ-mass}
In-plane slice images for the coupled lowres-LSQ-DTV alongside a profile plot demonstrating the effectiveness of
$f_2$ at extracting the image background. The dashed line in the $f_2$ image indicates the location of the profile plot data.}
\end{figure}

The data RMSE constraints, on the other hand, are more sensitive and these should be tuned on a case-by-case
basis. Here, these parameters are set by visual inspection, but for clinical application an automated procedure
must be devised to establish good settings. The settings for the three cases are listed in Table~\ref{tab:parms}.
As stated, these data RMSE constraints are the most critical parameters of coupled lowres-LSQ-DTV as they are the primary
controls on the image regularization. The looser, i.e. larger, these values are, the more regularized $f_1$ and $f_2$
become. The minimum values for $\epsilon_1$ and $\epsilon_2$ can be found by setting these parameters to zero,
then noting what data discrepancy is actually achieved when coupled lowres-LSQ-DTV is run to convergence. Backing off
from these minimum values regularizes their respective images $f_1$ and $f_2$. Visually, the $\epsilon$ values
are increased so that speckle noise is removed from $f_1$ and blobby artifacts on the scale of $d_2$ are removed
from $f_2$. The depth blur parameters involved in Eqs.~(\ref{disp1}) and (\ref{disp2}) depend on the features
of interest.

Finally, the PDHG algorithm parameters for coupled lowres-LSQ-DTV are set to the same values for each
of the three cases, as they also do not need to be finely tuned,
\begin{equation*}
\gamma=5, \;\;\; \rho=1.75, \; \text{and} \; \; \beta = 100;
\end{equation*}
these parameters are explained in Sec.~\ref{app:pdhg}.

\begin{figure}
\begin{center}
\begin{tabular}{c}
\includegraphics[width=0.9\columnwidth]{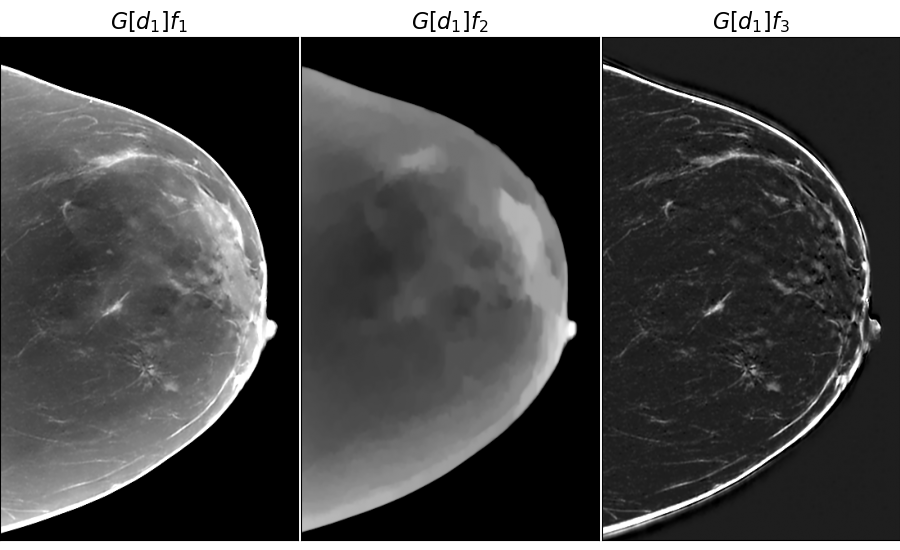}
\end{tabular}
\end{center}
\caption[Low-resolution coupled lowres-LSQ-DTV architectural distortion images.]{ \label{fig:DTV-LSQ-distortion}
In-plane slice images for the coupled lowres-LSQ-DTV for the patient case containing an architectural distortion.}
\end{figure}

\section{Results on wide-angle DBT image reconstruction for three patient cases}
\label{sec:results}

We demonstrate the application of coupled lowres-LSQ-DTV on three
patient cases acquired on a Siemens Mammomat B.brilliant scanner.
Approval of the Institutional Review Board (IRB) of Washington University (Protocol ID: 202410112) was obtained for this retrospective,
Health Insurance Portability and Accountability Act-compliant study conducted from data found in
institutional patient databases. Written informed consent was waived by the IRB. Three cases are selected that display
different features of clinical interest: a solid mass, an architectural distortion, and a calcification cluster.

\begin{figure}
\begin{center}
\begin{tabular}{c}
\includegraphics[width=0.95\columnwidth]{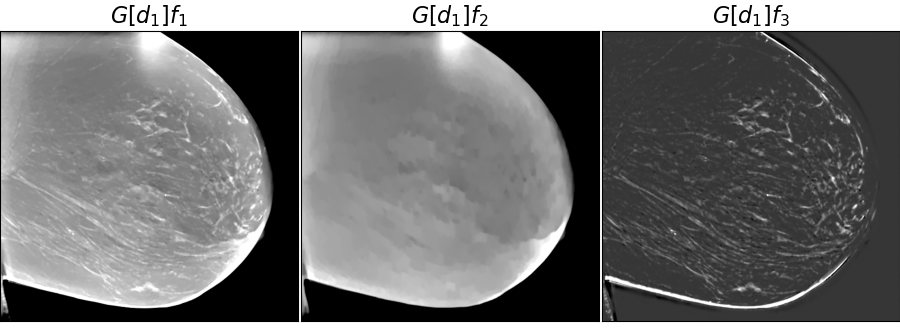}
\end{tabular}
\end{center}
\caption[Low-resolution coupled lowres-LSQ-DTV calcification cluster images.]{ \label{fig:DTV-LSQ-calcs}
In-plane slice images for the coupled lowres-LSQ-DTV for the patient case containing a calcification cluster.}
\end{figure} 

Results for coupled lowres-LSQ-DTV 
are shown in Fig.~\ref{fig:DTV-LSQ-mass}, where an in-plane slice image is shown that contains
a suspicious mass near the chest wall. The line profile and difference image show how $f_2$ is able to model the
 image background without following
clinically relevant features such as the mass, vessels, or fibroglandular tissue;
additionally, $f_2$ is able to account for background low-frequency artifacts.
Coupled lowres-LSQ-DTV is also applied
to a patient case with a distortion and the results are shown in Fig.~\ref{fig:DTV-LSQ-distortion}.
The feature of interest, the architectural distortion, is better visualized in the high-resolution region-of-interest (ROI)
images shown
in Fig.~\ref{fig:distortionROIs}, but Fig.~\ref{fig:DTV-LSQ-distortion} shows how the background $f_2$ image
picks up on the uniform dense region just inside the nipple, allowing for the difference image $f_3$ to highlight
fine structures in this region. Finally, the same algorithm is applied to a case with a significant microcalcification
cluster, and results are displayed in Fig.~\ref{fig:DTV-LSQ-calcs}. Of particular interest is that coupled
lowres-LSQ-DTV is able to account for the truncation artifacts at the top of the image, and consequently
these artifacts are greatly reduced in the $f_3$ difference image. The ROI focusing on the calcifications is shown
in Fig.~\ref{fig:calcsROIs}.

\begin{figure}
\begin{center}
\begin{tabular}{c}
\includegraphics[width=0.8\columnwidth]{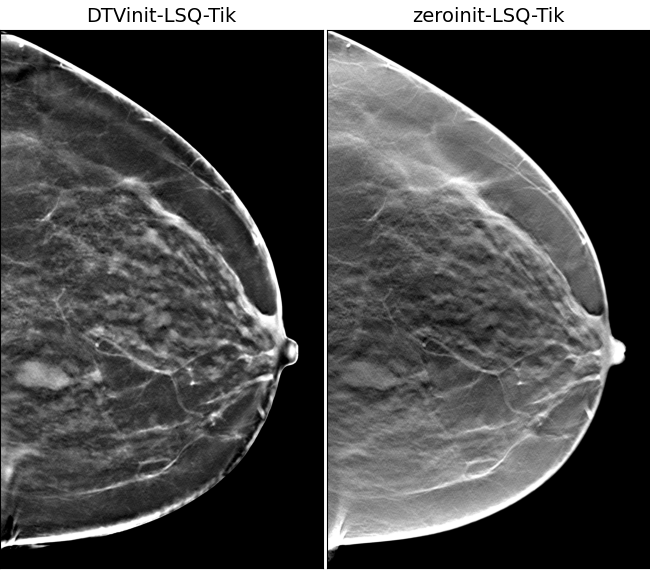}
\end{tabular}
\end{center}
\caption[High-resolution mass image comparison.]{ \label{fig:LSQ-Tik-mass}
In-plane images of DTVinit-LSQ-Tik (left) and zeroinit-LSQ-Tik (right) for the same slice shown in Fig.~\ref{fig:DTV-LSQ-mass}.}
\end{figure}

\begin{figure}
\begin{center}
\begin{tabular}{c}
\includegraphics[width=0.8\columnwidth]{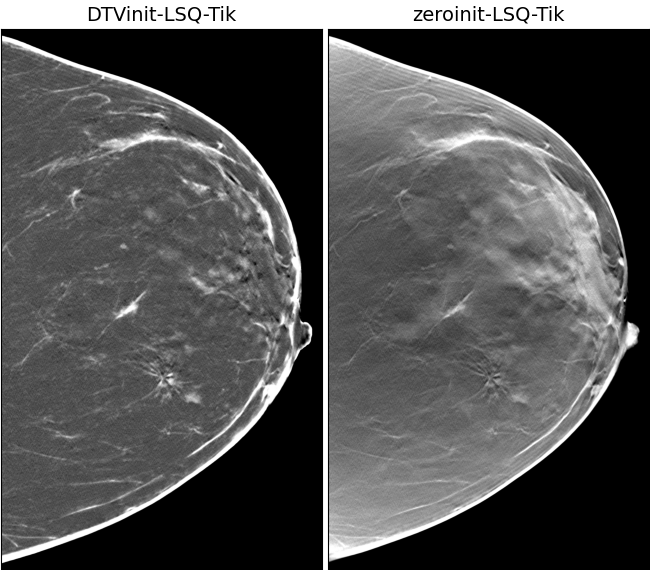}
\end{tabular}
\end{center}
\caption[High-resolution architectural distortion image comparison.]{ \label{fig:LSQ-Tik-distortion}
In-plane images of DTVinit-LSQ-Tik (left) and zeroinit-LSQ-Tik (right) for the same slice shown in Fig.~\ref{fig:DTV-LSQ-distortion}.}
\end{figure} 

In order to obtain standard DBT resolution, we employ DTVinit-LSQ-Tik as described in Sec.~\ref{sec:highres}.
The results of DTVinit-LSQ-Tik are compared with zeroinit-LSQ-Tik, which is found by the same algorithm except that
the initial image $h_0$ has all voxel values set to zero. No background artifact removal for zeroinit-LSQ-Tik is performed.
Example in-plane slice images are shown in Fig.~\ref{fig:LSQ-Tik-mass}. Because the background low-frequency artifacts are removed
from DTVinit-LSQ-Tik, it is possible to display the image in a tighter grayscale window than the image for zeroinit-LSQ-Tik.
Interestingly, the improved depth-resolution afforded by the coupled lowres-LSQ-DTV initial image resolves the tangle of glandular tissue
in the outer breast (image top) and improved in-plane contrast highlights the mass at the chest wall. 

\begin{figure}
\begin{center}
\begin{tabular}{c}
\includegraphics[width=0.9\columnwidth]{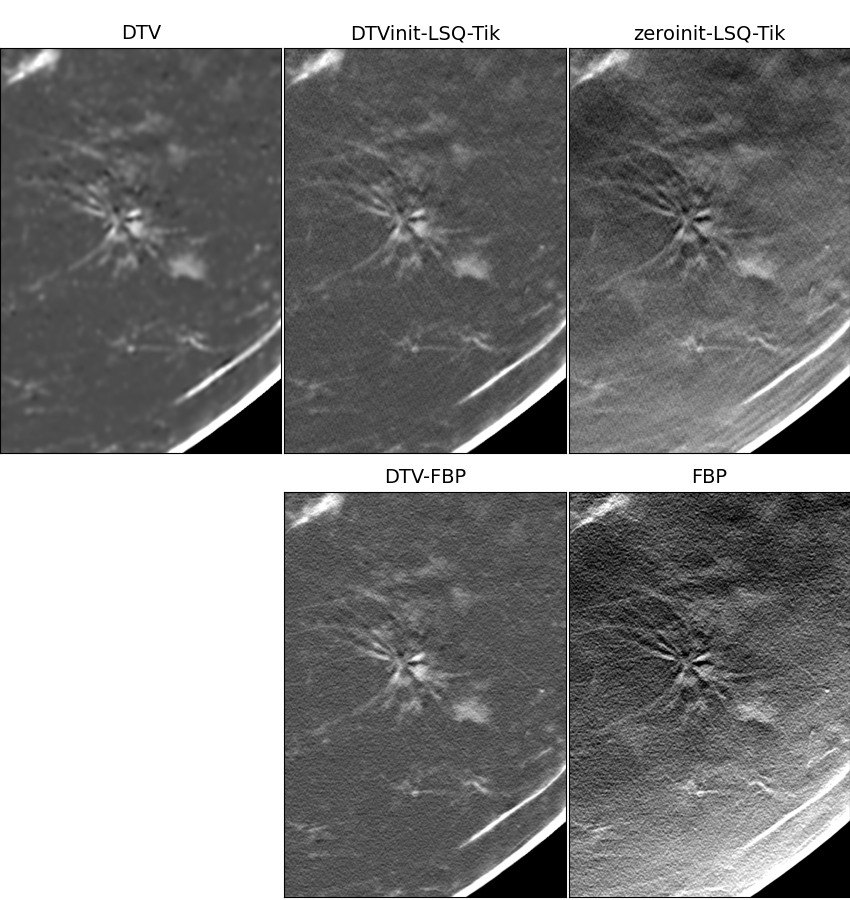}
\end{tabular}
\end{center}
\caption[High-resolution architectural distortion ROI image comparison.]{ \label{fig:distortionROIs}
Blow-up ROIs of the architectural distortion from the case displated in Fig.~\ref{fig:LSQ-Tik-distortion}.
The single image in the left column shows the input  coupled lowres-LSQ-DTV result. The middle column shows
the results from the proposed high-resolution algorithms. The right column shows the ROIs for the comparison
algorithms.}
\end{figure}

A comparison between DTVinit-LSQ-Tik and zeroinit-LSQ-Tik for the distortion case is shown in Fig.~\ref{fig:LSQ-Tik-distortion}.
Once again the the former algorithm is seen to remove much of the display slice clutter due to structure outside of the viewing
plane. The removal of low-frequency gray scale variation allows a tight display window thereby improving contrast of important
features of clinical interest such as the architectural distortion. A blow-up of this feature is shown in Fig.~\ref{fig:distortionROIs}
along with a different high-resolution image generated by high-resolution DTV-FBP. Recall that DTVinit-LSQ-Tik and DTV-FBP
images are generated by Eq.~(\ref{disp1}) and Eq.~(\ref{disp2}), respectively.
In this figure, the initial coupled lowres-LSQ-DTV ROI is shown so that the increase in resolution
by either high-resolution DTVinit-LSQ-Tik or DTV-FBP can be appreciated. Both of these results are compared against
their counterpart algorithms, zeroinit-LSQ-Tik and FBP respectively, where the advantage in distortion contrast
and reduced clutter can be clearly seen. The two alternative algorithms are presented because there is a clear
difference in noise texture and image sharpness at this level of magnification. 

\begin{figure}
\begin{center}
\begin{tabular}{c}
\includegraphics[width=0.95\columnwidth]{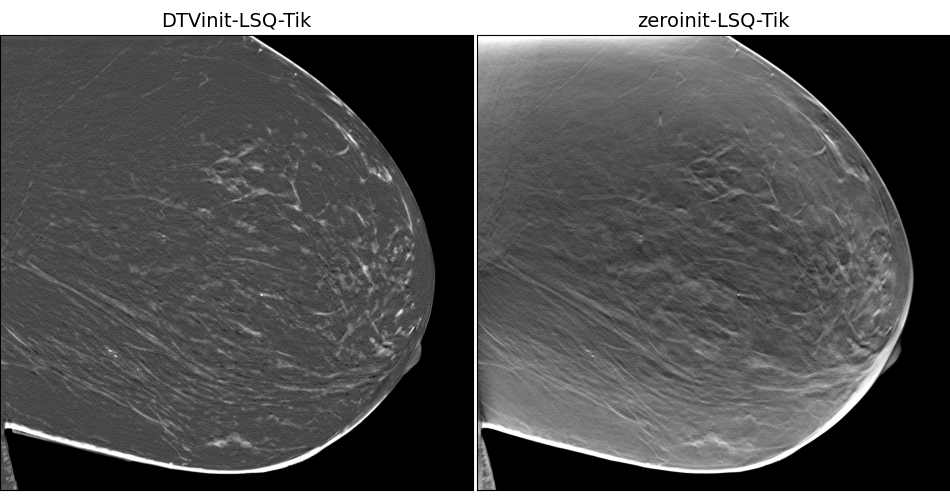}
\end{tabular}
\end{center}
\caption[High-resolution calcification cluster image comparison.] { \label{fig:LSQ-Tik-calcs}
In-plane images of DTVinit-LSQ-Tik (left) and zeroinit-LSQ-Tik (right) for the same slice shown in Fig.~\ref{fig:DTV-LSQ-calcs}.}
\end{figure}

High-resolution images for the case with micro-calcifications are shown in Fig.~\ref{fig:LSQ-Tik-calcs}.
As before there is less clutter in the DTVinit-LSQ-Tik image; the new aspect highlighted by this case is that
the truncation artifacts seen at the top left in the zeroinit-LSQ-Tik image is effectively suppressed in the DTVinit-LSQ-Tik
image. Focusing on the microcalcifications, ROI images for this case are shown in
Fig.~\ref{fig:calcsROIs}. The difference between the low-resolution ROI and the high-resolution images for
this ROI is greater than was the case for the distortion images, due the the fact that calcifications are smaller features.
As before, there is a difference in noise texture and sharpness between the high-resolution DTVinit-LSQ-Tik and DTV-FBP
images. Use of coupled lowres-LSQ-DTV, however, does not have a large impact for microcalcification imaging as there
is not a large difference between the middle and right columns in terms of micro-calcification visualization.

\begin{figure}
\begin{center}
\begin{tabular}{c}
\includegraphics[width=0.9\columnwidth]{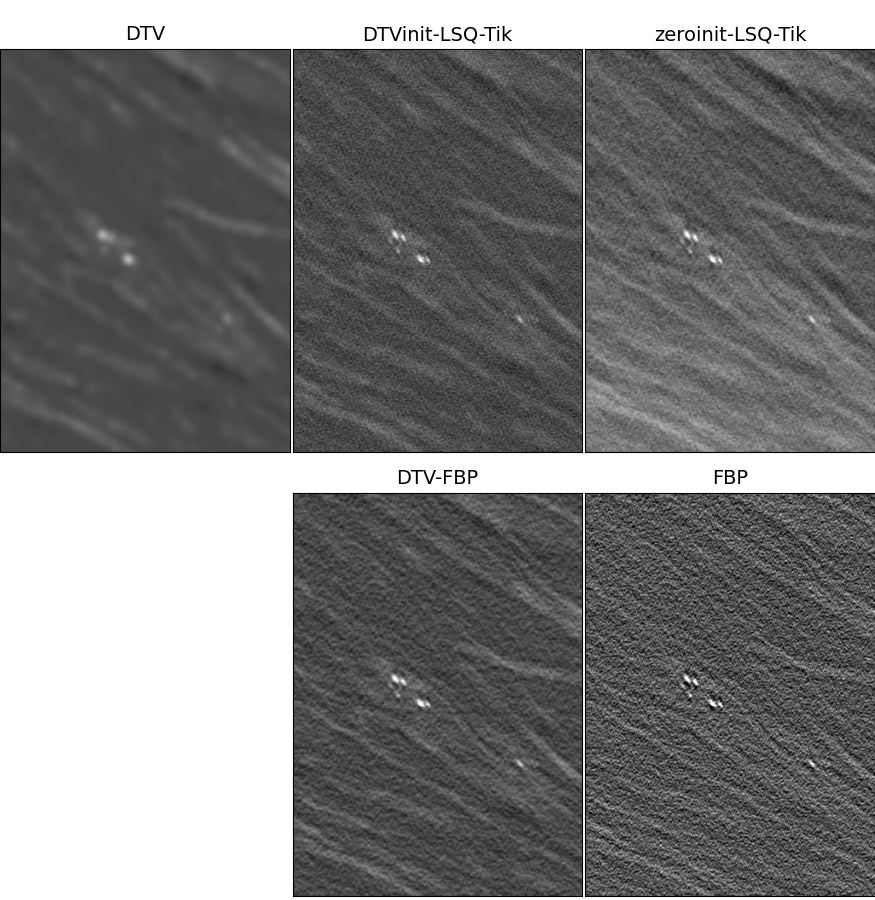}
\end{tabular}
\end{center}
\caption[High-resolution calcification cluster ROI image comparison.]{ \label{fig:calcsROIs}
Blow-up ROIs of a calcification cluster from the case displated in Fig.~\ref{fig:LSQ-Tik-calcs}.
The ROI panel presentation is the same as that of Fig.~\ref{fig:distortionROIs}.}
\end{figure} 

\section{Discussion}
\label{sec:discussion}
We have demonstrated, with three patient cases, the use of coupled lowres-LSQ-DTV in increasing
depth-resolution and slice uniformity for improving contrast and visualization of solid masses and spiculations indicating
an architectural distortion. Visualization of micro-calcifications is also shown, but little change in their
conspicuity is noted, because existing image reconstruction techniques can already resolve such features well.
The proposed two-step reconstruction: using accurately solved coupled lowres-LSQ-DTV followed by either high-resolution
FBP or LSQ-Tik, with truncated iteration, is highly efficient. In the shown examples, the low-resolution problem
is a factor of 512 smaller than the high-resolution problem, because the volume is represented on a grid with
$x$ and $y$ dimensions reduced by a factor of four each and the $z$ dimension by a factor of two.
Additionally, the projection dimensions $u$ and $v$ are
reduced by a factor of four each. The computationally demanding steps are X-ray projection and back-projection.
The low-resolution coupled lowres-LSQ-DTV algorithm executes a total of 4,000 of these projection operations for the
1,000 iterations used to achieve accurate solution; each iteration has two projection and back-projection
operations due to the simultaneous processing of $f_1$ and $f_2$ images.
 The high-resolution LSQ-Tik algorithm is run for 10 iterations,
involving 20 projection operations. But because of the
factor of 512 gain in efficiency from the reduced grid size, the coupled lowres-LSQ-DTV phase takes half of the computational effort
of high-resolution LSQ-Tik. Of course, high-resolution DTV-FBP requires much less computational than effort than LSQ-Tik because
it involves only a single back-projection operation.

While the proposed algorithms are seen to be effective in improving depth-resolution and in-plane contrast, there
are a few directions of investigation that would demonstrate, or further improve, clinical utility.
Of immediate concern is the bright skin-line; we have shown the direct output of the developed image reconstruction
algorithms without any post-processing. Skin-line correction is, however, an important post-processing step that
is needed for clinical use. Another point is to devise non-linear image reconstruction in the second high-resolution stage
specifically for denoising and improving micro-calcification conspicuity \cite{zheng2017detector}, noting that
this work by Zheng {\it et al.} was performed on a wide-angle system with a 60 degree scan arc.
Finally, a larger study of patient cases is needed where the imaging properties can be studied as a function of
breast density category.

\section{Conclusion}
\label{sec:conclusion}

We have taken the opportunity of this special issue honoring Dr. Barrett to report on the latest developments on our work
in digital breast tomosynthesis image reconstruction. The goal of the algorithm design is to improve depth resolution, in-plane contrast,
and to remove non-uniformity artifacts.
While the results from the shown clinical cases show promise, it is clear that
more studies are needed to demonstrate improved utility. Much of the future effort in this direction will involve carrying out
image quality studies promoted by Dr. Barrett; DBT simulations for
clinical tasks sensitive to image reconstruction algorithms and parameters,
e.g. detection of subtle signals such as microcalcifications
or spiculations on the margin of a suspicious mass, will be formulated
and employed to assess the algorithms in an objective manner, as we have done
in Refs.~\citenum{sanchez2016use,rose2017investigating}.

The contributions Dr. Barrett has made permeate modern image science.
The impact of his work
goes far beyond the usual metrics of academic productivity, i.e. number of citations or successful grant applications. His work
has been so foundational that it often goes uncited because it is part of the bedrock of the field.
Thus we are grateful for the opportunity to spell out, in this work on DBT, the many ways
that his contributions have influenced our work.

\section{Appendix: PDHG implementation of Eq.~(\ref{DTVcopt3D})}
\label{app:pdhg}

This appendix shows how to fit Eq.~(\ref{DTVcopt3D}) into the generic PDHG form
\begin{equation}
\label{genericOpt}
\min_x \left\{ F(K x) + G(x) \right\},
\end{equation}
where the $F$ and $G$ functions are convex and possibly non-smooth; $K$ is a linear transforms;
and $x$ is the unknown vector.
With a multiterm objective function such as the one presented in Eq.~(\ref{DTVcopt3D}) the function
$F$ ends up representing a sum of convex functions and $K$ is formed by stacking multiple matrix blocks,
which are the linear transforms inside of each of the objective terms. As shown in Sec.~2 of
Ref.~\citenum{sidky2026efficient}, we make the following associations
\begin{equation*}
x =
\left[
\begin{array}{c}
f_1 \\
f_2 \\
f_3
\end{array}
\right], \;\;\;
y = 
\left[
\begin{array}{c}
y_1 \\
y_2 \\
\vdots \\
y_{13}
\end{array} \right] = Kx, \; \; \;
\hat{F}(y)=  \sum_{i=1}^{13} \hat{F}_i(y_i) =  \sum_{i=1}^{13} F_i \left(y_i \|K_i\|_2 /\nu_i \right),
\end{equation*}
\begin{equation*}
K = \left[
\begin{array}{ccc}
\nu_1 K_1/ \|K_1\|_2 & \cdot & \cdot \\
\nu_2 K_2/ \|K_2\|_2  & \cdot & \cdot\\
\vdots  & \vdots & \vdots\\
\nu_6 K_6/ \|K_6\|_2 & \cdot & \cdot\\
\cdot & \nu_7 K_7/ \|K_7\|_2 & \cdot \\
\cdot & \nu_8 K_8/ \|K_8\|_2 & \cdot\\
\vdots  & \vdots & \vdots\\
\cdot & \nu_{12} K_{12}/ \|K_{12}\|_2 & \cdot\\
\cdot & \cdot & \nu_{13} K_{13}/ \|K_{13}\|_2 
\end{array} \right],
\end{equation*}
where $\| K_i \|_2$ is the matrix norm, i.e. largest singular value, of the matrix $K_i$ and $\nu_i$
are matrix scaling parameters that can be optimized for PDHG convergence efficiency.
The function and matrix block assignments are as follows
\begin{align*}
F_1(y_1) &=  \delta(y_1 \; | \; \|y_1 - R[c]g\|_2 \le \epsilon_1 \cdot \sqrt{N_\text{data}}), \;\;\;
K_1 = R[c]XG[d_1] \\
F_7(y_7) &=  \delta(y_7 \; | \; \|y_7 - R[c]G_\text{det}[d_d]g\|_2 \le \epsilon_2 \cdot \sqrt{N_\text{data}}), \;\;\;
K_7 = R[c]XG[d_2], \\
F_2(y_2) &= \alpha_x\|y_2\|_1, \;\;\; K_2 = \partial_x, \;\;\;F_8(y_8) = \alpha_x\|y_8\|_1, \;\;\; K_8 = \partial_x, \\
F_3(y_3) &= \alpha_y\|y_3\|_1, \;\;\; K_3 = \partial_y, \;\;\;F_9(y_9) = \alpha_y\|y_9\|_1, \;\;\; K_9 = \partial_y, \\
F_4(y_4) &= \alpha_a\|y_4\|_1, \;\;\; K_4 = \partial_a, \;\;\;F_{10}(y_{10}) = \alpha_a\|y_{10}\|_1, \;\;\; K_{10} = \partial_a, \\
F_5(y_5) &= \alpha_b\|y_5\|_1, \;\;\; K_5 = \partial_b, \;\;\;F_{11}(y_{11}) = \alpha_b\|y_{11}\|_1, \;\;\; K_{11} = \partial_b, \\
F_6(y_6) &= \alpha_1\|y_6\|_1, \;\;\; K_6 = I, \;\;\;F_{12}(y_{12}) = \alpha_1\|y_{12}\|_1, \;\;\; K_{12} = I, \\
F_{13}(y_13) &= \alpha_3\|y_{13}\|_1, \;\;\; K_{13} = I,\\
G(x) &= \delta([f_1,f_2,f_3] \; | \; f_1-f_2-f_3 = 0),
\end{align*}
where the indicator function $\delta( x \; | \; C )$ is 0 ($\infty$) if $x$ does (does not) satisfy the equality or inequality $C$;
the oblique directional derivatives denoted with $a$ and $b$ are a combination of partials in the $x$ and $z$ direction
\begin{equation*}
\partial_{\theta} = \cos \theta \partial_x + \sin \theta \partial_z,
\end{equation*}
where $a$ and $b$ denote the directions $\theta = -25^\circ$ and $25^\circ$, respectively.

\begin{algorithm}
\hrulefill
\begin{algorithmic}[1]
\State $x^{(k+1)} \leftarrow (I + \tau \partial G)^{-1} \left( x^{(k)} - \tau \left(\sum_i \nu_i \hat{K}_i^\top \lambda_i^{(k)}\right) \right) $
\State $\tilde{x} = 2 x^{(k+1)} - x^{(k)} $
\State $\lambda_i^{(k+1)} \leftarrow (I + \sigma_i \partial \hat{F}_i^\star)^{-1} \left( \lambda_i^{(k)} + \sigma_i \nu_i \hat{K}_i \tilde{x} \right)  \;\;\; \forall i$
\end{algorithmic}
\hrulefill
\caption{PDHG update steps for solving Eq.~(\ref{genericOpt}).
The normalized matrices $\hat{K}_i = K_i/\|K_i\|_2$ are used for conciseness.
}
\label{alg:pdhg}
\end{algorithm}

As explained in Ref.~\citenum{sidky2026efficient}, the algorithm for solving Eq.~(\ref{genericOpt}) involves three
update steps, which are repeated in Algorithm~\ref{alg:pdhg} for completeness. Lines 1 and 3 of this algorithm involve
possibly unfamiliar notation which is fully explained in Ref.~\citenum{sidky2026efficient}; the reason why the update
steps are given is to specify the complete set of step-size parameters. Line 1 of Algorithm~\ref{alg:pdhg} involves
the scalar step-size $\tau$ and the matrix scaling parameters $\nu_i$; Line 3 involves the scalar step-sizes $\sigma_i$.
These step-size parameters are chosen so that
\begin{equation}
\label{stepcondition}
\tau \left\| \sum_i \sigma_i \nu^2_i \hat{K}^\top_i \hat{K}_i \right \|_2 = 1,
\end{equation}
which is based on the convergence condition for PDHG. Because the index $i$ runs through all 13 matrix blocks, there are
a total of 27 parameters ($\tau$, $\sigma_i$, and $\nu_i$) to tune for optimizing the PDHG convergence rate. The condition
in Eq.~(\ref{stepcondition}) reduces this number to 26 independent parameters, but tuning 26 parameters is still a daunting task.

In Ref.~\citenum{sidky2026efficient}, a framework for X-ray tomography is proposed that simplifies the parameter selection.
In that work, it is proposed to provide heuristic arguments for selecting step-size weight parameters, $\hat{w}_i$, and
ratio parameters $r_i$, which are related to $\tau$, $\sigma_i$, and $\nu_i$ by the equations
\begin{equation*}
w_i = w \hat{w}_i = \tau \sigma_i \nu_i, \; \text{and} \; \; r_i = \sigma_i/\tau,
\end{equation*}
where $\hat{w}_i$ are relative weights and the scalar magnitude $w$ is determined by Eq.~(\ref{stepcondition})
\begin{equation}
\label{wmag}
w = \left\| \sum_i \hat{w}_i \hat{K}^\top_i \hat{K}_i \right\|^{-1}_2.
\end{equation}
In Ref.~\citenum{sidky2026efficient}, it is proposed to set all step-size ratios equal to a single scalar
$\beta$ and to weight all linear transforms $K_i$ the same except for the matrix blocks involving the X-ray transform,
which should have greater weighting $\gamma$. For Eq.~(\ref{DTVcopt3D}) there are two matrix blocks, $K_1$
and $K_7$, that involve the X-ray transform. Accordingly, we make the following assignments
\begin{align*}
\hat{w}_1 & = \hat{w}_7 = \gamma, \\
\hat{w}_i & = 1, \; \; \forall i \neq 1\text{ or }7, \\
r_i & = \beta,
\end{align*}
where $\gamma$ is selected in the interval $\gamma \in [1, \infty]$, and $\beta$ is a positive real number.
As a result, the search space for the step-size parameters is reduced from 26 dimensions down to 2, which
is much more manageable. As is done in Ref.~\citenum{sidky2026efficient}, the reduction in step-size parameter search
space allows the use of an accelerated predictor-corrector PDHG algorithm that comes with an additional
extrapolation parameters $\rho$.
The particular instatiation of the PDHG algorithm calls for derivation of proximal operators involving
the function terms $F_i(y_i)$ and $G(x)$ and the details on how to derive those steps are presented
in Ref.~\citenum{sidky2026efficient}.

\subsection*{Disclosures}
The authors have no relevant financial interests in the manuscript and no other potential conflicts of interest to disclose.

\subsection* {Code, Data, and Materials Availability} 
There is no code or data in this work for sharing.

\subsection* {Acknowledgments}
This work is supported in part by NIH Grant R01-CA287302 and the
RSNA Research and Education Foundation (RR24-342).
The contents of this article
are solely the responsibility of the authors and do not necessarily represent the official views of the
National Institutes of Health.






\listoffigures

\end{spacing}
\end{document}